
%
\input amstex
\mag=\magstep1
\documentstyle{amsppt}
\NoRunningHeads
\vsize=23.5true cm
\hsize=15.5true cm
\topmatter
\title Crepant Resolution of Trihedral Singularities
\endtitle
\author Yukari ITO \endauthor
\affil Department of Mathematical Sciences \\
 University of Tokyo \\
 7-3-1 Hongo, Bunkyo, Tokyo, 113, Japan \\
 yukari-i{\@}tansei.cc.u-tokyo.ac.jp \endaffil
\endtopmatter
\document
\subhead  \nofrills{\bf
\S 1. Introduction}
\endsubhead \par
The purpose of this paper is to  construct a crepant resolution of quotient
singularities by finite subgroups of $SL(3,\Bbb C)$ of certain type, and prove
that each Euler number of the minimal model is equal to the number of conjugacy
classes. \par
 The problem of finding a nice resolution of quotient singularities by finite
subgroups of $SL(3,\Bbb C)$ arose from mathematical physics. In the superstring
theory, the dimension of the space-time is 10, four of them are usual space and
time dimensions, and other six are compactified on a compact Calabi-Yau space
$M$.  From a point of view of algebraic geometry, the Calabi-Yau space is a
smooth three-dimensional complex projective variety whose canonical bundle is
trivial and
fundamental group is finite.\par
 In the physics of superstring theory, one considers the string propagation on
a
manifold $M$ which is a quotient by a finite subgroup of symmetries $G$. By a
physical argument of string vacua of $M/G$, one concludes that the correct
Euler number for the theory should be
the ``orbifold Euler characteristic"[3,4], defined by \par
$$ \chi(M,G) = \frac 1{\vert G \vert} \sum_{gh=hg} \chi(M^{<g,h>}),$$
 where the summation runs over all pairs of commuting elements of $G$,
and $M^{<g,h>}$ denotes the common fixed set of $g$ and $h$. For the physcist's
interest, we only consider $M$ whose quotient space $M/G$ has trivial canonical
  bundle.

\proclaim{ Conjecture I} {\rm([3,4])}\par
\it{There exists a resolution of singularities $\widetilde{M/G}$ s.t.
$\omega_{\widetilde{M/G}} \simeq \Cal O_{\widetilde{M/G}}$, and
$$ \chi(\widetilde{M/G}) = \chi(M,G). $$} \endproclaim
This conjecture follows from its local form [5]:

\proclaim{ Conjecture II} {\rm(local form)}\par
\it{Let $G\subset SL(3,\Bbb C)$ be a finite group. Then there exists a
resolution
 of singularities\par $\sigma$ : $\widetilde X \longrightarrow {\Bbb C}^3/G$
with
$\omega_{\widetilde X} \simeq \Cal O_{\widetilde X} $ and
$$ \chi(\widetilde X) = \sharp \{ \text{conjugacy class of } G \}. $$
} \endproclaim
In algebraic geometry, the conjecture says that a minimal model of the quotient
space by a finite subgroup of $SL(3,\Bbb C)$ is non-singular.\par

Conjecture II  was proved for abelian groups  by Roan ([13]), and
independently by Markushevich, Olshanetsky and Perelomov ([8]) by using toric
method. It was also proved for 5 other groups, for which $X$ are hypersurfaces:
(i) $W{A_3}^+, W{B_3}^+, W{C_3}^+,$ where $WX^+$ denotes the positive
determinant part of the Weyl group $WX$ of a root system $X$, by Bertin and
Markushevich
 ([1]), (ii) $H_{168}$  by Markushevich ([7]), and
(iii) $I_{60}$  by Roan ([14]).

 \par
In this paper, we prove Conjecture II for solvable groups of type
 $(C)$ .
\definition{ Definition} {\it Trihedral group} is a finite group
 $G=<H,T>\ \subset SL(3,\Bbb C)$,
where $H \subset SL(3,\Bbb C)$ is a finite group generated by diagonal matrices
and
$$ T=\pmatrix
     0 & 1 & 0 \\
     0 & 0 & 1 \\
     1 & 0 & 0 \\
    \endpmatrix $$
\enddefinition
\definition{Definition}
{\it Trihedral singularities} are quotient singularities by trihedral groups.
\enddefinition
\definition{Definition}
A resolution of singularities $f:Y \longrightarrow X$ of a normal variety $X$
is {\it crepant} if $K_Y=f^*K_X$.
\enddefinition
\proclaim{Theorem 1.1 \rm{( Main Theorem )}}\par
\it{Let $X=\Bbb C^3/G$ be a quotient space by a trihedral group $G$. Then there
exists a crepant resolution of singularities
$$f:\widetilde X \longrightarrow X,$$
 and
$$ \chi(\widetilde X) = \sharp \{ \text{conjugacy class of }G \}. $$}
\endproclaim

Trihedral singularities are 3-dimensional version of $D_n$-singularities, and
they are non-isolated and many of them are not complete intersections.
Their resolutions are similar to those of $D_n$-singularities.
There is a nice combination of the toric resolution and Calabi-Yau
resolution.\par

By the way, the conjecture II  is true in dimension 2 (i.e., the case of
$SL(2,\Bbb C)$ (Cf.[5])), but in the case of $SL(4,\Bbb C)$ there exists a
counter-example; in the
case of group $G= <[-1,-1,-1,-1]>$ (diagonal matrix), which is a finite
subgroup
of $SL(4,\Bbb C)$, but there isn't a crepant resolution.
\par
We describe the contents of this paper.
In section 2, we recall the results on 3-dimensional canonical singularities,
and see some properties among the invariants under the group action.\
In the next three sections, we prove the Main Theorem. In section 3, we find a
normal subgroup $G'$ of $G$ and calculate the order of $G'$.
In section 4, we construct a crepant resolution of trihedral singularities with
data on $G'$.
In section 5, we compute the Euler number of $\widetilde X$ and see that it
coincides with the number of the conjugacy classes.
The list of finite subgroups of $SL(3,\Bbb C)$ and some examples are in
UTMS-94-18 (preprint of University of Tokyo).\par

 The author would like to express her hearty gratitude to Professor Y.~Kawamata
for his valuable advice and encouragement, and to Professors M.~Reid,
N.~Nakayama and  M.~Kobayashi for their helpful discussions and encouragement.
She would also like to thank Professors K-i.~Watanabe and T.~Uzawa for useful
advices concerning the list of subgroups of $SL(3,\Bbb C)$.
\subhead \nofrills{\bf
\S 2. Some properties of quotient singularities by $SL(n,\Bbb C)$
}\endsubhead
\par
\definition{Definition}([11]) A normal variety $X$ has {\it canonical}
(resp.{\it terminal}) {\it singularities} if it satisfies the following two
conditions:\par
\roster
 \item"({\bf i})" for some integer $r\ge 1$,
                the Weil divisor $rK_X$ is Cartier;
 \item"({\bf ii})" if $f : Y \longrightarrow X$ is a resolution of $X$ and
${E_i}$ the exceptional prime divisors of $f$, then
        $$ rK_X = f^*(rK_X) + \sum a_iE_i\  \text{  with }\  a_i \ge 0\
(\text{resp. } a_i > 0)$$
\endroster
\enddefinition

\definition{Definition}
In the formula $K_Y=f^*K_X+\sum a_iE_i, a_i $ is called the {\it discrepancy}
at $E_i$, and the minimum number of $r$ is called the {\it index}.
\enddefinition
\example{Remark 2.1} Canonical singularities of index 1 are the same as
rational Gorenstein singularities.\endexample

\proclaim{Proposition 2.2} \par
\it{Let $G$ be a finite subgroup $\not= \{1\}$ of $SL(3,\Bbb C)$, and $X$ be a
quotient space of
${\Bbb C}^3$ by $G$. Then $X$ is canonical but not terminal } \endproclaim
\demo{\bf{Proof}}
It follows from the classification of 3-dimensional terminal singularities
([10], [9]) and the
classification of invariant rings of finite subgroups of $SL(3,\Bbb C)$ which
are complete intersections (Cf.[15]). \qed
\enddemo
We obtain this result without the classifications, if $X$ is a hypersurface:
\proclaim{Proposition 2.3} \par
\it{Let $X$ be a hypersurface defined by $(f=0)$ in $\Bbb C^n$,
$(x_1, \cdots, x_n)$ a coordinate of $\Bbb C^3$,
$\alpha(x_i):=$ a  weight of $x_i$ and $\alpha(f):= \text{min}\{ \alpha(x^m) |
x^m \ \text{is a monomial in} \ f\}$.\par
If $$ \sum_{i=1}^n \alpha (x_i) = \alpha(f) + 1, $$
then $X$ is not terminal.\par
Furthermore if $X$ is a quotient space by $G \subset SL(3,\Bbb C)$, then it is
canonical but not terminal.}
\endproclaim
\demo {\bf{Proof}}
This follows from the fact ([12]) that
$$ X:\text{terminal (resp. canonical)}\ \Longrightarrow \sum_{i=1}^n\alpha
(x_i)\ > \ \alpha(f) +1 \ \text{(resp. $\geq$)}  $$
\enddemo

 In the following sections $(\S 3, \S 4, \S 5)$, we prove Theorem 1.1.

\subhead \nofrills {\bf \S 3. Normal subgroup $G'$} \endsubhead \par
Let $G'$ be the subgroup of the group $G=<H,T>$ consisting of all the diagonal
matrices. Then $G'$ is a normal subgroup, and an abelian group. We consider the
order of $G'$.
\proclaim{Proposition 3.1} \par
\it{ $|G'|$ is one of the following holds. \roster
\item $|G'| \equiv 0$ (mod 3)
\item $|G'| \equiv 1$ (mod 3)
\endroster} \endproclaim

\demo{bf{Proof}}
 Let $\frac 1{r}(a,b,c),$  $a+b+c \equiv 0$ (mod $r$) be a type of an element
$g=[\theta^a,\theta^b,\theta^c]\in G'$ where $\theta$=exp$(2\pi \sqrt{-1}/r)$.
 Then there exist three possibilities of the type of the element of $G'$:
\roster\item $a\not=b,  b\not=c,  c\not=a$
\item $a=b,b\not=c$
\item $a=b=c$
\endroster
The number of the elements of the types (1) and (2) is divisible by 3, because
there exists three elements
$$ \frac 1{r}(a,b,c), \frac 1{r}(b,c,a), \frac 1{r}(c,a,b). $$
So it is sufficient to consider the number of the elements of type (3).
It is at most three, because they are
$$ \frac 1{3}(1,1,1), \frac 1
{3}(2,2,2),\text{ and identity element}$$
So, if $G'$ has an element $\frac 1{3}(1,1,1)$, then the order of $G'$ is
congruent to zero modulo 3. Otherwise, $|G'|$ is congruent to 1 modulo 3.
\qed \enddemo
{}From now, we call the type of $G'$ as the following:\par
Type (I)  when $|G'| \equiv 1$ (mod 3) \par
Type (II)  when $|G'| \equiv 0$ (mod 3)
\pagebreak
\subhead \nofrills {\bf \S 4. Crepant resolution of trihedral singularities }
\endsubhead \par

 \proclaim {Proposition 4.1} \par
\it{Let $X = {\Bbb C} ^3/ G$, and $Y = {\Bbb C}^3 /G'$. Then there exists the
following diagram: }
$$ \CD
   @. @. {\widetilde X} \\
   @. @. @VV{\tau}V     \\
   @. {\widetilde Y} @>{\widetilde\mu}>> {\widetilde Y}/{\frak A_3} \\
   @. @VV{\pi}V @VV{\widetilde\pi}V \\
   {\Bbb C^3} @>>> {\Bbb C^3}/{G'}=Y @>{\mu}>> {\Bbb C^3}/G=X
 \endCD $$
\it{where $\pi$ is a resolution of the singularity of $Y$, and $\widetilde\pi$
is the induced morphism, $\tau$ is a resolution of  the singularity by $\frak
A_3$, and $\tau \circ \widetilde \pi$ is a crepnt resolution of the singularity
of $X$.}\endproclaim

\demo{\bf{Sketch of the proof}} As a resolution $\pi$ of $Y$, we take a toric
resolution, which is also crepant. Then we lift the $\frak A_3$-action on $Y$
to itts minimal resolution $\widetilde Y$ and form the quotient $\widetilde
Y/\frak A_3$. This quotient gives in a natural way a partial resolution of the
singularities of $X$. The minimal resolution $\widetilde X \longrightarrow
\widetilde Y/\frak A_3$ of the singularities
of $\widetilde Y/\frak A_3$ induces a complete resolution of $X$.\par
Under the action of  $\frak A_3$, the singularities of $\widetilde Y/\frak A_3$
lie in the union of the image of the exceptional diviser of $\widetilde Y$
under $\widetilde Y  \longrightarrow \widetilde Y/\frak A_3$ and the locus
$C:(x=y=z)$. \par
In the resolution $\widetilde Y$ of $Y$, the group $\frak A_3$ permutes
exceptional divisors. So the fixed points on the exceptional divisors consist
of one point or three points. \enddemo

Now, we see the proof more precisely.
\par
At first, we recall a toric resolution when $G'$ is a finite abelian group in
$SL(3,\Bbb C)$.\par
Let $\Bbb R^3$ be a vector space, and $\{ e^i | i=1,2,3\}$ the standard base.
For all $v=\frac 1{r}(a,b,c)\in G'$, $L= $the lattice generated by $e^1,e^2$
and $e^3$,  $N:=L\oplus \sum \Bbb Zv$,
$$ \sigma = \left\{ \sum_{i=1}^3 x_ie^i \in \Bbb R^3, \quad x_i \geq 0, \forall
i \right\}. $$
We regard $\sigma$ as a rational convex polyhedral cones in $N_{\Bbb R}$. The
corresponding affine torus embedding $X_{\sigma}$ is defined as Spec$\Bbb
C[\check{\sigma} \cap M]$, where $M$ is the dual lattice of $N$ and
$\check{\sigma}$ is a dual cone of $\sigma$ in $M_{\Bbb R}$ defined by
$\check{\sigma} =\{ \xi\in M_{\Bbb R} | \xi(x) \geq 0 , \forall x \in \sigma\}
$. $X_{\sigma}$ is isomorphic to $\Bbb C^3/G'$ as a algebraic varieties. \par
$\qquad \qquad \varDelta :=$ the simplex in $N_{\Bbb R}$
$$ \left\{\sum_{i=1}^3 x_ie^i  \quad ; x_i \geq 0,\quad \sum_{i=1}^3 x_i=1
\right\} \qquad \qquad$$ \par
$$  t:N_{\Bbb R} \longrightarrow \Bbb R \qquad \sum_{i=1}^3x_ie^i \longmapsto
\sum_{i=1}^3 x_i \ \qquad \qquad \qquad  $$\par
$$ \Phi := \left\{ \frac 1{r}(a,b,c)\in G'\  | \ a+b+c=r \right\}
\qquad \qquad \qquad $$
\proclaim{Lemma 4.2} \par
\it{$Y=\Bbb C^3/G'$ corresponds to $\sigma$ in $N=L\oplus \sum_{v \in \Phi}\Bbb
Zv$.}
\endproclaim
\demo{\bf{Proof}}
Since $Y=\text{Spec}(\Bbb C[x,y,z]^{G'})$, $x^iy^jz^k$ is $G'$-invariant if and
only if $ \alpha i+\beta j + \gamma k \in \Bbb Z$ for all
$(\alpha,\beta,\gamma) \in G'$.  \qed \enddemo
\example{Remark 4.3}
Let $\frac 1{r}(a,b,c) \not= (0,0,0)$ be an element of $G'$. There are two
types.\roster
\item \quad $abc \not= 0$
\item \quad $abc = 0$
\endroster
We denote by $G_1$ (resp. $G_2$) the set of the elements of type (1) (resp.
(2)). So $G'\backslash \{e\} = G_1 \amalg G_2$. \par
There are also two types in $\Phi$ as above. So we denote by $\Phi_1$ (resp.
$\Phi_2)$ the set of lattice points of type (1) (resp. (2)). Then $\Phi=\Phi_1
\amalg \Phi_2$.\par
Let $\lambda_i$ be maps from $G_i$ to $\Phi_i$ (i=1,2).
$$ \lambda_i\  :\ G_i \ \longrightarrow \ \Phi_i $$
where $\lambda_1$ maps $g=\frac 1{r}(a,b,c)\ (a+b+c=r)$ and $g^{-1}= \frac
1{r}(r-a,r-b,r-c)$ to a lattice point $\frac 1{r}(a,b,c)$, and $\lambda_2$ maps
$g=\frac 1{r}(a,b,c)$ to a lattice point $\frac 1{r}(a,b,c)$.
 $$ \left\{G_1\right\} \overset {2:1} \to \longrightarrow
\left\{\Phi_1\right\}\ , \qquad \left\{G_2\right\} \overset {1:1} \to
\longrightarrow \left\{\Phi_2\right\} $$
 \par
Therefore there exist two correspondence between the sets of elements in $G'$
and the set of exceptional divisors.
\endexample
\proclaim{Cliam I} \it{There exists a toric resolution of $Y$ where $\frak A_3$
acts symmetrically on the exceptional divisors.}
\endproclaim
\demo{\bf{Proof}}  We show that we can construct a simplicial decomposition
$\{ \sigma(s)\}_{s \in S}$ of the simplex $\varDelta$ which is $\frak
A_3$-invariant and the set of  vertices is exactly $\Phi \cup \{e^i\}_{i=1}^3$.
\par
 Let us consider the distance $d$ between $\frac 1{r}(a,b,c)$ and $\frac
1{3}(1,1,1)$ given by
$$ d\left(\frac 1{r}(a,b,c)\right)=\left\vert\frac a{r}-\frac 1{3}\right\vert +
\left\vert\frac b{r}-\frac 1{3} \right\vert
                        + \left\vert \frac c{r}- \frac 1{3}\right\vert.$$
Case of type (I)
\roster
\item Find a lattice point $P= \frac 1{r}(a,b,c)$ whose distance $d$ is the
minimum among the points in $\Phi$ in the domain $D=\{x,y\geq z$\}.
\item Make a triangle whose vertices are $P$,\ $P'= \frac 1{r}(b,c,a)$  and
$P''= \frac 1{r}(c,a,b)$.
\item There is no lattice point inside or on the boudary of the triangle
$PP'P''$ because of the minimarity of $d$.
\item Decompose a quadrangle $Q_1$ whose vertices are $P$,\ $P'$,\ $(1,0,0)$
and  $(0,1,0)$, into  simplexes using vertices in $\Phi$.
 We call this decomposition $S_1$.
\item By the action of $\frak A_3$, we obtain $S_2$ and $S_3$ on 2 other
quadangles $Q_2,\ Q_3$. Therefore we obtain a ``symmetric"
resolution.
\endroster \par
Case of type (II)
\roster
\item Take a lattice point $P$ from $D\backslash \{ \frac 1{3}(1,1,1) \}$ as
(1) in above case.
\item Make a triangle $PP'P''$ as above, connect between the vertices and a
point $\frac 1{3}(1,1,1)$ with lines.
\item Similar as (4),(5) in above case.
\endroster
Therefore we obtain the result.
\qed \enddemo
\pagebreak

\proclaim{Claim II}
\it{Let $X_S$ be the corresponding torus embedding, then $X_S$ is non-singular.
}\endproclaim
\demo{\bf{Proof}}
It is sufficiant to show that the $\sigma(s)$ are basic.
Let $ w^1,w^2,w^3 \in \Phi \cup \{e^i\}_{i=1}^3$ which are linearly independent
over $\Bbb R$. Assume that the simplex
$$ \left \{ \sum_{i=1}^3 \alpha_iw^i \ | \ \alpha_i \geq 0 ,\ \sum_{i=1}^3
\alpha_i=1 \right \} $$
intersects $\Phi \cup \{e^i\}_{i=1}^3$ only at $\{w^i\}_{i=1}^3.$
\par
Lattice $N_0$ generated by $\{w^i\}_{i=1}^3$ is sublattice of $N$.
If we assume $N\not= N_0$, then there exists
$\beta=\beta_1w^1+\beta_2w^2+\beta_3w^3 \ \in\  N\backslash N_0$\  $( 0 <
\beta_i  < 1,\  \beta_i \in \Bbb R).$\par
$t(\beta)=\sum\beta_i t(w^i) = \sum \beta_i$, $0 < \sum \beta_i < 3$ and $t(N)
\in \Bbb Z$, then $t(\beta)=1$ or 2.
If $t(\beta)=2$, then we can replace it by $\beta'=
\sum_{i=1}^3(1-\beta_i)w^i$,  so we can assume that $t(\beta)=1$. \par
Now, there exists an element  $\beta$ in $
\{ \sum \alpha_iw^i \ | \ \alpha_i \geq 0, \ \sum \alpha_i =1\} \cap (N -
N_0)$, which is contained in $\varDelta \cap N$.
Since $$N = \left\{\bigcup \Sb v \in \Phi \endSb ( v\oplus L)\right\}\bigcup L,
$$
$\varDelta \cap N = \Phi \cup \{e^i\}_{i=1}^3.$ From our assumption, we
conclude that $\beta=w^i$ for some $i$, which contradicts $\beta\not \in N_0$.
Therefore $N=N_0$.
\qed \enddemo
We obtain a crepant resolution $\pi_S:X_S\longrightarrow \Bbb C^3/G'.$

\proclaim{Claim III}
\it{Let $F$ be a fixed locus on $\widetilde Y$ under the action of $\frak A_3$,
and $E$ be exceptional divisors of $\widetilde Y \longrightarrow Y$. Then
$$ F_0 := F \cap E = \cases \text{1 point \qquad \ if $G'$ is type (I)} \\
                                 \text{3 points \qquad if $G'$ is type (II)}.
\endcases $$}
\endproclaim
\demo {\bf{Proof}}
Considering the dual graph of exceptional divisors by the toric resolution and
from Remark 4.3, we can identify the three exceptional divisors except the
central component. Then there are 2 possibilities of the central component;
\par
 Type (I): one point.\par
 Type (II): a divisor which is isomorphic to $\Bbb P^2$.
Let $(x:y:z)$ be a coordinate of $\Bbb P^2$, then the action of $\frak A_3$ is
$$ (x:y:z) \longmapsto (y:z:x). $$
Then the number of the fixed points are three, whose coordinates are
$$(1:1:1)\ ,\quad (1:\omega:\omega^2)\ , \quad (1:\omega^2:\omega).\qed$$
\enddemo
Futhermore the $\frak A_3$-action in the neighbourhood of a fixed point is
analytically isomorphic to some linear action. \par
Now, we consider the resolution of the singularity of $\Bbb C^3$ by the group
$\frak A_3$.\par
$F_0$ consists of 1 or 3 points, and $F$ = $F_0 \cup C'$ where $C'$ is a strict
transform of the fixed locus $C$ in $Y$ under the action of $\frak A_3$. \par
Under the action
$$ T :(x,y,z) \longmapsto (y,z,x),$$
the fixed locus $C$ is defined by $x=y=z$, then one of the fixed points in
$F_0$ lies in the strict transform $C'$, and  the singular locus in $\widetilde
Y /\frak A_3$ are
$$\cases \text{ $G'$ is  type (I)}\  \quad \Longrightarrow \quad C' \\
         \text{ $G'$ is type (II)} \quad \Longrightarrow \quad C'\cup \{\text{2
points} \}
\endcases$$
\proclaim {Claim IV} \it{Let $Z=\Bbb C^3/\frak A_3$, then $\chi(\widetilde Z) =
\chi( \Bbb C^3, \frak A_3) = 3.$}
\endproclaim
\demo{\bf{Proof}}
There are two representations of $\frak A_3$ in $SL(3,\Bbb C)$:
$$ \text{(i)}\quad T' = \pmatrix
         1 & 0 & 0 \\
         0 & \omega & 0 \\
         0 & 0 & \omega^2 \\
        \endpmatrix
 \qquad \quad
    \text{(ii)} \quad T'' = \pmatrix
         \omega & 0 & 0 \\
         0 & \omega & 0 \\
         0 & 0 & \omega \\
         \endpmatrix     $$
The quotient singularities by $T'$ are $A_2 \times \{x-$axis\} which are not
isoalated and the exceptional divisors
of the resolution are two $\Bbb P^1$-bundle intersecting at their section. \par
On the other hand, the quotient singularity by $T''$ is isolated and the
exceptional divisor of the resolution is isomorphic to $\Bbb P^2$. Therefore
each Euler  number is 3.
\qed \enddemo

If $G'$ is of type (I), $\widetilde Y/\frak A_3$ has a non-isolated singularity
of Type (i). If $G'$ is of type (II), $\widetilde Y/\frak A_3$ has one
non-isolated singularity of type (i) and two isolated singularities of type
(ii).

\proclaim {Claim V} \it{
The resolution  $\tau \circ \widetilde \pi $ is a crepant resolution.}
\endproclaim

\subhead \nofrills{\bf \S 5. Correspondence with conjugacy class} \endsubhead
\par
\proclaim {Lemma 5.1} \par
\it{Let $X:= \Bbb C^3/<G',T>$, and $f:\widetilde X \longrightarrow X$ the
crepant resolution as above. Then the Euler number of $\widetilde X$ is given
by
$$ \chi(\widetilde X) = \frac 1{3} (\vert G' \vert -k)+3k $$
where $$ k  =
       \cases
           1 \qquad \text{if\ \   $|G'| \equiv 1$ (mod 3) \  (type (I))}  \\
           3 \qquad \text{if\ \   $|G'| \equiv 0$ (mod 3) \  (type (II))}
       \endcases $$}
\endproclaim

\demo{\bf{Proof}}
For an abelian group $G'$, we have a toric resolution
$$ \pi : \widetilde Y \longrightarrow Y = \Bbb C^3/G', $$
and $\chi(\widetilde Y) = \vert G'\vert$ ([8],[11]). \par
By the action of $\frak A_3$, the number of fixed points in the exceptional
divisor by $\sigma$ is equal to $k$, hence
$$ \chi(\widetilde Y/\frak A_3)=\frac 1{3} (\vert G' \vert -k)+k $$
By the resolution of the fixed loci, Eular characteristic of the each
exceptional locus is 3. (Claim IV) \par
Therefore,
$$ \chi(\widetilde X) = \frac 1{3}(\vert G'\vert -k) + 3k  \qed$$
\enddemo
Now, we compute the number of the conjugacy class of the group $G$.
We observe the order of the normal subgroup $G'$. \par

\proclaim{Theorem 5.2} \par
\it{$$\chi (\widetilde X) = \sharp \{ \text{conjugacy class of $G$} \} $$
}\endproclaim
\demo{\bf{Proof}}
\roster
 \item  Case (I) : $\vert G'\vert = 3m+1 , (0 <m\in \Bbb Z)$
\par For the nontrivial element $g \in G'$,
there are three conjugate elements
$g$,$T^{-1}gT$,$TgT^{-1}$.
There are $m$ triples of this type.
There are 3 other  conjugacy classes $[e],[T],[T^{-1}]$.
 Therefore, there are $m+3$ conjugacy classes in $G$. \par
Then
$$ \allowdisplaybreaks \align
   \chi(\widetilde X) &= \frac 1{3}(|G'|-1)+3 \\
                      &= m+3 \\
                      &= \sharp \{ \text{ conjugacy class of } G \}
\endalign $$
 \item  Case (II) : $\vert G'\vert = 3m , (0 <m \in \Bbb Z)$
\par There are 3 elements in the center of $G'$: $e, a=[\omega,\omega,\omega],
a^2=[\omega^2,\omega^2,\omega^2]$. The other  $3m-3$ elements in $G'$ are
divived into $m-1$ conjugacy classes as in (1). There are 6 other conjugacy
classes
$[T],[T^{-1}],[aT],[a^2T],[aT^{-1}]$ and $[a^2T^{-1}]$.
Therefore, there are $m+8$ conjugacy classes in $G$. \par
Then
$$ \allowdisplaybreaks \align
     \chi(\widetilde X) & =\frac 1{3} (\vert G'\vert -3) +9 \\
                        & = m+8 \\
                        & = \sharp \{ \text{ conjugacy class of } G \} \qed
\endalign     $$
\endroster
\enddemo
Therefore, Main theorem (Theorem 1.1) is proved!

\enddocument